\documentclass{article}
\usepackage[utf8]{inputenc}
\usepackage{amsmath}
\usepackage{xcolor}
\usepackage[
backend=biber,
style=alphabetic,
sorting=ynt
]{biblatex}
\usepackage[affil-it]{authblk}
\usepackage{graphicx}
\usepackage{subcaption}
\usepackage{textcomp}
\providecommand{\keywords}[1]
{
  \small	
  \textbf{\textit{Keywords---}} #1
}

\title{Exploring the Dynamics of Virulent and Avirulent Aphids: A Case for a ``Within Plant" Refuge}
\author[1]{Aniket Banerjee%
  \thanks{Electronic address: \texttt{aniketb@iastate.edu}; Corresponding author}}
\affil[1]{Department of Mathematics, Iowa State University, Ames, IA 50011, USA}
\author[2]{Ivair Valmorbida%
}
\affil[2]{Department of Entomology, Iowa State University, Ames, IA 50011, USA}
\author[2]{Matthew O'Neal%
}
\author[1]{ Rana Parshad%
}
\date{May 2021}

\begin{document}

\maketitle
\begin{abstract}
    The soybean aphid, \textit{Aphis glycines} (Hemiptera: Aphididae), is an invasive pest that can cause severe yield loss to soybeans in the northcentral United States. A tactic to counter this pest is the use of aphid-resistant soybean varieties. However, the occurrence of virulent biotypes can alter plant physiology and impair the use of this management strategy. Soybean aphids can alter soybean physiology primarily by two mechanisms, feeding facilitation and the obviation of resistance, favoring subsequent colonization by additional conspecifics. We developed a non-local, differential equation population model, to explore the dynamics of these biological mechanisms on soybean plants co-infested with virulent and avirulent aphids. We then use demographic parameters from laboratory experiments to perform numerical simulations via the model. These simulations successfully mimic various aphid dynamics observed in the field. Our model showed an increase in colonization of virulent aphids increases the likelihood that aphid-resistance is suppressed, subsequently increasing the survival of avirulent aphids, producing an indirect, positive interaction between the biotypes. These results suggest the potential for a “within plant” refuge that could contribute to the sustainable use of aphid resistant soybeans.
\end{abstract}

\keywords{soybean, biotypes, IRM, non-local ODE model.}

\section{Introduction}
The challenge in using insect-resistant plants is the development of sub-populations of the target pest that survive on these plants. Efforts to limit their frequency can be achieved through the use of insecticide resistant management (IRM) plans (Tabashnik et al. 2013). One subset of pests that are particularly challenging to manage with resistant plants are aphids (Hemiptera: Aphididae). This is due, in part, to their capacity to evolve virulent biotypes that can survive on resistant plants and their remarkable capacity to reproduce asexually across multiple generations (Crowder and Carriere, 2009).  The ability for some aphids to manipulate their host plant so that sub-populations survive (Elzinga and Jander 2013, Varenhorst et al. 2015), leaves open the possibility that current models for explaining aphid-plant interactions may not accurately describe unique features of some systems, like that of the soybean aphid, \textit{Aphis glycines} (Hemiptera: Aphididae) on soybeans.\\
The soybean aphid was first detected in 2000, and has become one of the most important insect pests of soybean in the major production areas of the Midwest US (Tilmon et al. 2011). It has a heteroecious holocyclic life cycle that utilizes a primary host plant for overwintering and a secondary host plant during the summer. In the spring, aphids emerge and produce three asexual generations on common buckthorn, \textit{Rhamnus cathartica} L., then migrate to soybeans \textit{Glycine max} L. Aphids continue to reproduce asexually on soybeans, producing as many as 15 generations during the summer (Ragsdale et al. 2011). In North America, aphids arrive on soybean fields in June, where populations increase by four orders of magnitude and at the end of the growing season (mid-September), aphids begin the migration back to their overwintering host, reproduce sexually and overwinter in the egg stage (Ragsdale et al. 2004). Within Iowa, populations of aphids large enough to reduce soybean yield occurred on 40\% of growing seasons from 2004 to 2019, with populations peaking in the middle to the end of August (Dean et. al. 2020).\\
Colonization and feeding by an insect herbivore can alter the plant’s physiology, favoring the subsequent colonization of additional conspecifics (Price et al. 2011). There are two mechanisms by which this susceptibility can be induced, feeding facilitation and the obviation of resistance. Feeding facilitation is a more general mechanism by which the general physiology of the host plant is altered by the herbivore, often in a density dependent manner. A more specific mechanism that inducts susceptibility is the obviation of traits that confer resistance to the herbivore (e.g., Baluch et. al. 2012, Varenhorst et al. 2015). This mechanism requires a subset of the herbivores population that is virulent, capable of surviving on the resistant genotype of the host plant. By obviating the resistance through a physiological change to the plant, avirulent subpopulations can now survive on the resistant plant. Both mechanisms allow sub-populations that vary by genotype (i.e., virulent and avirulent) to co-exist on resistant host plants. These mechanisms have been observed in populations of soybean aphids, when colonizing soybean plants that have resistance to soybean aphids (O’Neal et al. 2018, Varenhorst et al. 2015). Field surveys in North America have observed that soybean aphid biotypes can co-occur in the same fields (Cooper et al. 2015, Alt et al. 2019). Laboratory studies have also shown that virulent and avirulent biotypes can co-exist on a shared plant for at least 2-3 generations (Varenhorst et al. 2015). However, there is no empirical evidence that soybean aphid biotypes can co-occur on the same plant throughout a growing season.\\
Exploring the potential for aphid biotypes to co-exist on a shared plant is challenging with empirical studies, as several abiotic and biotic factors can affect aphid populations dynamics. We explored if avirulent and virulent soybean aphids can co-exist on a shared soybean plant by expanding upon existing population models developed for aphids. Included within this expanded model is resistance to aphids in the host plant, and features that reduce this resistance consistent with feeding facilitation and obviation of aphid-resistance. We calculated population growth rates for both virulent and avirulent aphids on susceptible and resistant plants, and use these parameters within this model. The model was used to determine if a 'within-plant' refuge is possible, such that avirulent aphids can persist on a resistant plant throughout a growing season. Finally, we discuss how the addition of a third trophic layer (i.e., predators) may affect the outcome of this model and its implications for the implementation of IRM.

\section{Material and Methods}
\subsection{Model development from a single biotype}
Kindlemann et al. 2010 are among the first to propose a model for the population dynamics of aphids using a set of differential equations (Model 1).

\begin{eqnarray}
\label{eq:m1}
  \frac{dh}{dt}&=&ax ;  h(0)=0 \nonumber\\ 
     \frac{dx}{dt}&=&(r-h)x \ ;  x(0)=x_0\nonumber\\
     \nonumber
\end{eqnarray}

Where $h(t)$ is the cumulative population density of a single aphid biotype at time t; $x(t)$ is the population density at time $t$, $a$ is a scalar constant, and $r$ is the growth rate of the aphids. The aphid population initially rises due to the linear growth term, but as the cumulative density becomes greater than the growth rate $r$, the population starts to decrease, due to the effects of competition. This results in a hump shaped population density over time, typical of a boom-bust type scenario (Kot, 2001). This is an apt description of aphid dynamics, particularly when exploring soybean aphid dynamics on soybeans during the growing season. The type of population growth described by this model (Kindlemann et al. 2010) has been observed in soybean aphids in North America (e.g., Brosius et al. 2007, Catangui et al. 2009), with colonization in June, then a gradual buildup of population, peaking in August and declining with aphids dispersing in September to their overwintering host.\\
The model (1) is quite different from the classical logistic growth model, which predicts growth to a certain carrying capacity. It is an example of a non-autonomous model, wherein the right-hand side of the differential equation depends explicitly on time. The rigorous mathematical analysis of such systems is quite involved, and the methods of classical autonomous systems do not apply (Langa et al. 2002). Hence the rigorous dynamical analysis of model (1) is not found elsewhere. However, it provides a starting point to model more intricate aphid dynamics, particularly when a species presents two or more biotypes.

\subsection{Virulent and avirulent aphids: Two biotypes}

Sub-populations of a herbivorous species can be organized into biotypes, defined as genotypes capable of surviving and reproducing on host plants containing traits (e.g., antibiosis and or antixenosis) conferring resistance to that herbivore (Downie, 2010). Specifically, for the soybean aphid, biotypes are classified based on their ability to colonize soybean varieties expressing \textit{Rag}-genes (\textit{Rag} is derived from the expression, resistance to \textit{Aphis glycines}). For example, soybean aphid biotype 1 is susceptible to all \textit{Rag}-genes, therefore it is called avirulent. Biotype 2 is virulent to \textit{Rag1} (Kim et al. 2008), biotype 3 is virulent to \textit{Rag2} (Hill et al. 2010), and biotype 4 is virulent to both \textit{Rag1} and \textit{Rag2}, capable of surviving on plants with these genes either alone and together (Alt and Ryan-Mahmutagic, 2013). These four soybean aphid biotypes have been found throughout the soybean producing areas of the Midwest US (Cooper et al. 2015, Alt et al., 2019).\\
The model proposed by Kindlemann et al. (2010) cannot describe aphid population dynamics on a soybean plant colonized by both virulent and avirulent aphids. This is because that model does not account for competition or cooperation between the two biotypes. First, the virulent and avirulent are in direct competition for space, similar to interspecies competition. The virulent aphids are also in competition for space with other virulent aphids, as avirulent aphids are in competition for space with other avirulent aphids, a case of intraspecies competition. These are interaction can produce competition through direct effects. Furthermore, on a resistant plant both the avirulent and virulent aphids are able to weaken the plant’s defenses via feeding facilitation. However, for the avirulent aphid this only occurs if it arrives in sufficiently large numbers (Varenhost et al. 2015). Thus, there is a definite resistant level in the plant that is dependent on the initial density of colonizing avirulent aphids. If the avirulent aphids arrive in sufficient numbers above this level, they could colonize a resistant plant. On the other hand, the virulent biotype alters the plant by obviating the resistance (O’Neal et al. 2018), allowing both virulent and avirulent aphids to survive on \textit{Rag1+2} plants. The suppression of the plant’s resistance level by the virulent biotype, eases the colonization process for the avirulent biotype, which is an indirect form of cooperation at play. Thus, the plant’s resistance is a dynamic process, dependent on the presence and densities of these biotypes.\\
In model (2), our goal is to simulate a system that considers two soybean aphid biotypes (virulent and avirulent) that are attempting to colonize a soybean plant containing aphid resistance in the form of \textit{Rag} genes. The model includes all of the earlier mentioned interactions and considers a dynamic aphid-resistance level in the plant that is affected by the density of the aphids, in the capacities mentioned earlier.\\
The expanded model (2) is as follows:\\
\begin{eqnarray}
  \frac{dh}{dt}&=&a(x_A+ x_V ) \nonumber\\ 
    \frac{dx_A}{dt}&=&(r-h)(x_A-R) \nonumber\\ 
  \frac{dx_V}{dt}&=&(r-h)x_V  \nonumber\\
   \frac{dR}{dt}&=&-(k_r x_V+k_f x_V+k_f sgn(x_A-A) x_A)R \nonumber  \nonumber 
\end{eqnarray}

Here $x_A(t)$ refers to avirulent aphid population density and $x_V(t)$ refers to the virulent aphid population density, $h$ is the combined cumulative population density of both avirulent and virulent aphids, respectively, at time $t$. $r$ is the maximum potential growth rate of the aphids. $a$ is a scaling constant relating aphid cumulative density to its own dynamics. $R$ is the dynamic resistance threshold of the plant (Berec, 2004). This decreases due to both avirulent and virulent aphid density, that is $x_V$ and $x_A$. It is measured in the same units as aphid density. $k_f$ is the rate of feeding facilitation and $k_r$ is the rate of obviation of resistance (Varenhorst et al. 2015). Feeding facilitation occurs for avirulent aphids, if their population is above a threshold level $A$, Below this threshold, the effect of feeding facilitation by an avirulent population is negligible.  Previous studies have demonstrated that obviation of resistance is much more effective in suppressing the resistance than feeding facilitation (Varenhorst et al. 2015). Therefore, $k_r> k_f$ whenever both the effects take place simultaneously. $Sgn (xA-A)$ is a Boolean function returning 0 or 1. It returns 1 if the input is strictly positive (i.e. $(xA – A) > 0$) or else it returns 0. This function regulates whether avirulent aphids have enough initial population density to induce effect of feeding facilitation on the plant.\\

\subsection{Model Parameters}
We explored a time series analysis of the soybean aphid population dynamics on a soybean plant containing aphid resistance in the form of \textit{Rag1+2} genes. We used values in the model based on our understanding of the dynamics between the two biotypes. One of the most important parameters of our model is the growth rate of the aphids ($r$). This determines the timing of the boom-bust scenario along with the cumulative population density. The growth rate $r$ of the biotypes on resistant (\textit{Rag1+2}) and susceptible plant was estimated using a life table analysis. Treatments consisted of two factors, soybean cultivar (susceptible and \textit{Rag1+2}) and aphid biotypes (avirulent and virulent). Soybean aphids used in this experiment have been kept at Iowa State University, reproducing parthenogenically on soybeans. Aphids were kept in separated growth chambers under controlled conditions [25 ± 2 °C, 70\% RH and 16:8 (L:D)]. The avirulent aphid was reared on a susceptible soybean (LD14‐8007), while the virulent was reared on a \textit{Rag1+2} soybean variety (LD14‐8001). Soybean seeds were sown in 8-cm-diamter plastic plots using a soil mixture (Sungro Horticulture Products, SS\#1-F1P, Agawam, MA, USA). Plants were kept in a greenhouse [25 ± 5 °C and a photoperiod of 16:8 (L:D)], watered three times per week and fertilized weekly after emergence (Peters Excel Multi-Purpose Fertilizer, 21-5-20 NPK).\\
Twenty-four hours before the beginning of the experiment, a single mix-aged, apterous adult aphid was transferred onto a soybean leaflet kept in a Petri dish within a growth room [25 ± 2 °C, 50\% RH and 16:8 (L:D)]. This allowed us to synchronize the age of the aphids ($\leq$24 h old). A total of 40 aphids were used for each avirulent and virulent soybean aphid. After 24 h, a single first instar nymph was transferred to the first trifoliate leave of a V2 (Fehr et al. 1971) soybean plant. Plants were then covered with a mesh net to prevent aphids from escaping and moving to another plant, were kept in a growth room [25 ± 2 °C, 50\% RH and 16:8 (L:D)], and watered three times per week. Each treatment combination consisted of 25 potted plants.\\
Evaluations were performed daily until the aphid died. Morphological characteristics were used to determine growth stage (Zhang, 1988, Voegtlin et al. 2004) and exuviae were removed once detected. When the aphids became adults, their offspring were counted and removed daily. All the aphids used in the treatment combination of avirulent aphid and \textit{Rag1+2} variety died within three days and were not included in the statistical analysis. Biological and demographic parameters were calculated using the TWOSEX-MSChart (Chi, 2020) program following the age-stage, two-sex life table theory (Chi and Liu, 1985; Chi, 1988). Means and standard error of population parameters were estimated using a bootstrap procedure (Huang et al., 2013) with 100 000 replicates. Differences among treatments were analyzed using a paired bootstrap test at 5\% significant level using the TWOSEX-MSChart program.\\

\section{Results}
\subsection{Parameters used in the model}
We adjusted the parameters within our model to explore several scenarios in which both virulent and avirulent aphids colonized a single plant. These scenarios are outlined in Table 1. The biological and demographic parameters of avirulent and virulent soybean aphids are presented in Table 2.  The treatment consisting of avirulent aphids on a  resistant plant was not include in the analysis because all the avirulent aphids died within three days. The finite rate of increase ($\lambda$=0.27) generated from the life table analysis was used as the growth parameter ($r$). The value used in the model is an average of the three treatments. In our time series analysis, the initial dynamic resistance ($R (0)$) and threshold avirulent population ($A$) is fixed as 30 aphids.  The scaling parameter $a$=0.000005 used in the model comes from Kindlmann et al. 2010. The variables that defined feeding facilitation ($k_f$) and obviation of resistance ($k_r$) were adjusted to explore the impact of each on the occurrence of the two biotypes.\\

\subsection{Feeding facilitation}
While studying feeding facilitation within the model, we removed the effect of obviation of resistance due to virulent aphids by setting $k_r=0$. The model shows that when virulent aphids are absent and the initial avirulent population is below the resistance threshold, the avirulent aphid goes to extinction as feeding facilitation does not take place due to a low initial population (Fig. 1A). However, when the population of avirulent aphid is above the resistance threshold, colonization takes place even in the absence of the virulent aphid (Fig. 1B). When the initial population of virulent aphids is very small and the population of avirulent aphids is below the resistance threshold, the effect feeding facilitation is insufficient to allow for the coexistence of both biotypes (Fig. 1C and 1D). Even in the absence of resistance obviation, if the initial population of virulent aphids is increased, the model accounts for this increase in the overall population of aphids on the plant by increasing the strength of feeding facilitation. This further sustains avirulent aphids on the plant throughout the season even when the initial avirulent aphid is below the resistance threshold (Fig. 1E and 1G).\\
A significant increase in the peak population of avirulent aphid is observed with an increase in their initial population. This relationship is apparent when comparing figures 1E and 1F, as an increase in the initial population of the avirulent aphid results in a higher maximum population. The model fixes the value of the host plant as a resource for aphids. Thus, an increase in the avirulent population results in the decrease in the maximum population attained by virulent aphids in the season.

\subsection{Obviation of resistance}

Obviation of resistance is a phenomenon by which the virulent aphids suppress the resistance of the resistant plant, allowing both the virulent and avirulent aphids to colonize and grow on a shared plant. In figure 2, we explored if our model could produce results consistent with this phenomenon. We set the initial population of the avirulent aphid lower than the resistance threshold in all the cases, while holding the initial avirulent population and the resistance threshold constant in all cases. This prevents feeding facilitation, as noted in the previous section.\\
As shown in figures 2A and 2B, when the virulent aphid population is very low (i.e., below $A$), the resistance is not yet suppressed and the avirulent aphid goes extinct. By increasing the initial virulent aphid population from 35 to 50, the avirulent aphid persists on the plant (Fig. 2C and 2D). There was a sufficient density of virulent aphids in these scenarios to suppress the resistance in the plant, allowing the avirulent aphid to survive on the resistant plant. As we increase the initial population of virulent aphid (Fig. 2E and 2F), the resistance declines much faster allowing the avirulent aphid to reach higher densities.\\
Obviation of resistance as described in this model, allows for the persistence of the avirulent aphid at population levels lower than what the resistance threshold would allow. With an increase in the initial virulent aphid, obviation of resistance begins sooner, resulting in higher populations of avirulent aphid across the modeled season. We also see that the closer the initial virulent aphid population is to the resistance threshold, the fewer initial virulent aphids are required to obviate resistance such that the avirulent aphid population is sustained.\\

\subsection{Coexistence of two biotypes}
If the initial avirulent population is higher or equal to the resistance threshold then the two biotypes of the aphids co-exist in the ecosystem for the whole season (Fig. 3). When initial population of avirulent aphid is set higher than the resistance threshold (Fig. 3) both biotypes can persist on the plant regardless of initial virulent aphid density. We also observe that if initial avirulent population is low but close to the resistance threshold then both populations co-exist (Fig. 2C-2F). Which population will dominate is dependent on the initial density of each biotype. For example, an increase of the virulent aphid's initial population results in a larger population density compared with the avirulent aphid population. Nonetheless, the peak for the population of both biotypes occurs at the same time.

\section{Discussion}
The model accounts for the dynamics of feeding facilitation with both avirulent and virulent soybean aphids. If the plant is resistant (i.e., it contains Rag-genes), then virulent aphids survive and feeding facilitation takes place at any population level. If this resistant plant is colonized by avirulent aphids only, survival is not guaranteed. However, if population of avirulent aphids is higher than the resistance threshold (R), feeding facilitation (Varenhorst et al. 2015) allows this biotype to survive on a resistant plant.\\
Currently, the frequency of virulent biotypes within North America is lower than that of avirulent biotypes (Cooper et al. 2015, Alt et al. 2019).  This scenario is necessary for \textit{Rag}-based resistance to remain useful for the management of this pest.  With the increased use of a \textit{Rag}-resistance, the frequency of virulent biotypes is expected to increase. Efforts to prevent virulent biotypes from increasing with increasing use of aphid-resistance is the goal of an IRM program. A strategy for the sustainable use of a resistant variety is the creation of a 'refuge' of susceptible plants that maintain a sufficient population of avirulent biotypes to reduce the frequency virulent biotypes in subsequent generations (Crowder and Carriere 2009).  By coupling a refuge with resistant plants that are so toxic that the virulence becomes functionally recessive, this strategy has preserved the value of insect-resistant plants (Tabashnik et al. 2013).  A challenge to this strategy is farmer acceptance of a practice that requires the cultivation of susceptible plants.  Refuges have been incorporated to the units of seed sold to farmers so that this strategy is practiced with limited input from the farmer. This practice has been referred to as a 'refuge-in-a-bag'. The soybean aphid/\textit{Rag}-resistance system suggests that a refuge could occur within a plant.  For this refuge to contribute to the management of virulent populations in an IRM program, the avirulent populations must persist throughout the growing season and contribute to the population that returns to the overwintering host.  This model suggests that there is the potential to contribute to such a refuge from the limited empirical evidence of feeding facilitation and the obviation or resistance. The full potential of a 'refuge-in-a-plant' requires more empirical evidence that the avirulent populations generated from a within plant refuge contribute to the overwintering population.\\
An additional detail that would improve our modeling efforts is the role of natural enemies in affecting the frequency of virulence. Natural enemies can affect the frequency of biotypes that are virulent to resistant host plants (Gould et al. 1991). A community of aphidophagous predators can be found in soybean fields that feed on soybean aphids during the summer. Studies manipulating aphid exposure to natural enemies using cages demonstrated that predators play a role in suppressing the growth of aphid populations in North America (Costamagna and Landis 2006, Costamagna et al. 2007, Bannerman et al. 2018). Parasitoids are an additional source of soybean aphid mortality that can significantly impact population growth (Frewin et al. 2010, Kaser and Heimpel 2018), especially in their native range (Liu et al. 2004). It is not known if the impact of either predators or parasitoids varies significantly by aphid biotype.  Revealing the impact of natural enemies on the dynamic relationship between biotypes on aphid-resistant plants may be critical for developing an IRM strategy for soybean aphids.\\

\section{Acknowledgements}
Funding was provided from an Insect Management Knowledge Program through Monsanto for the project titled “Extending durability of aphid resistance by understanding the mechanisms of virulence to Rag-genes in soybeans”. We also received support from the Agricultural Experiment Station at Iowa State University. AB and RP acknowledge valuable partial support via the National Science Foundation.

\section{References}

Alt, J., and M. Ryan‐Mahmutagic. 2013. Soybean aphid biotype 4 identified. Crop Sci. 53: 1491-1495.\\
Alt, J., M. Ryan, and D.W. Onstad. 2019. Geographic distribution and intrabiotypic variability of four soybean aphid biotypes. Crop Sci. 59: 84-91.\\
Baluch, S.D., H.W. Ohm, J.T. Shukle, and C.E. Williams. 2012. Obviation of wheat resistance to the Hessian fly through systemic induced susceptibility. J. Econ. Entomol. 105: 642-650.\\
Bannerman, J.A., B.P. McCornack, D.W. Ragsdale, N. Koper, and A.C. Costamagna. 2018. Predators and alate immigration influence the season-long dynamics of soybean aphid (Hemiptera: Aphididae). Biol. Control. 117: 87-98.\\
Berec, L.E. 2004. Mathematical modeling in ecology and epidemiology. Habilitation thesis, Masaryk University, Czech Republic.\\
Brosius, R.T., L.G. Higley, and T.E. Hunt. 2007. Population dynamics of soybean aphid and biotic mortality at the edge of its range. J. Econ. Entomol. 100: 1268-1275.\\
Catangui, M.A., E.A. Beckendorf, and W.E. Riedell. 2009. Soybean aphid population dynamics, soybean yield loss, and development of stage-specific economic injury levels. Agron. J. 101: 1080-1092.\\
Chi, H. 1988. Life-table analysis incorporating both sexes and variable development rates among individuals. Environ. Entomol. 17: 26-34.\\
Chi, H. 2020.TWOSEX-MSChart: A computer program for the age-stage, Twosex life table analysis. http://140.120.197.173/ecology/products.htm.\\
Chi, H., and H. Liu, H. 1985. Two new methods for the study of insect population ecology. Bull. Inst. Zool. Acad. Sin 24: 225-240. \\
Cooper, S.G., V. Concibido, R. Estes, D. Hunt, G.L. Jiang, C. Krupke, B. McCornack, R. Mian, M. O'Neal, V. Poysa, D. Prischmann‐Voldseth, D. Ragsdale, N. Tinsley, and D. Wang. 2015. Geographic distribution of soybean aphid biotypes in the United States and Canada during 2008–2010. Crop Sci. 55: 2598-2608.\\
Costamagna, A.C., and D.A. Landis. 2006. Predators exert top‐down control of soybean aphid across a gradient of agricultural management systems. Ecol. Appl.  16: 1619-1628.\\
Costamagna, A.C., D.A.  Landis, and C.D. Difonzo. 2007. Suppression of soybean aphid by generalist predators results in a trophic cascade in soybeans. Ecol. Appl. 17: 441-451.\\
Crowder, D.W., and Y. Carriére. 2009. Comparing the refuge strategy for managing the evolution of insect resistance under different reproductive strategies. Theor. Biol. 261: 423-430.\\
Dean, A.N., J.B. Niemi, J.C. Tyndall, E.W. Hodgson, and M.E. O'Neal. 2020. Developing a decision‐making framework for insect pest management: a case study using Aphis glycines (Hemiptera: Aphididae). Pest Manag. Sci. 77: 886-894.\\
Downie, D.A. 2010. Baubles, bangles, and biotypes: a critical review of the use and abuse of the biotype concept. J. Insect Sci.  10: 176.\\
Elzinga, D. and G. Jander. 2013. The role of protein effectors in plant-aphid interactions. Curr. Opinion Plant Biol. 16: 451-456. \\
Fehr, W.R. C.E. Caviness, D.T. Burmood, and J.S. Pennington. 1971. Stage of development descriptions for soybeans, Glycine max (L.) Merrill. Crop Sci. 11: 929–931.\\
Frewin, A.J., Y. Xue, J.A. Welsman, B.A. Broadbent, A.W. Schaafsma, and R.H. Hallett. 2010. Development and parasitism by Aphelinus certus (Hymenoptera: Aphelinidae), a parasitoid of Aphis glycines (Hemiptera: Aphididae). Environ. Entomol. 39: 1570-1578.\\
Gould, F., G.C. Kennedy, and M.T. Johnson. 1991. Effects of natural enemies on the rate of herbivore adaptation to resistant host plants. Entomolo. Exp. Appl. 58: 1-14.\\
Hill, C.B., L. Crull, T.K. Herman, D.J. Voegtlin, and G.L. Hartman. 2010. A new soybean aphid (Hemiptera: Aphididae) biotype identified. J. Econ. Entomol. 103: 509-515.\\
Huang, Y.B. and H. Chi. 2013. Life tables of Bactrocera cucurbitae (Diptera: Tephritidae): with an invalidation of the jackknife technique. J. Appl. Entomol. 137: 327-339.\\
Johnson, S.N., R.C. Rowe, and C.R. Hall.  2020. Aphid feeding induces phytohormonal cross-talk without affecting silicon defense against subsequent chewing herbivores. Plants. 9:1009.\\
Kaser, J.M., and G.E. Heimpel. 2018. Impact of the parasitoid Aphelinus certus on soybean aphid populations. Biol. Control. 127: 17-24.\\
Kim, K.S., C.B. Hill, G.L. Hartman, M.R. Mian, and B.W. Diers. 2008. Discovery of soybean aphid biotypes. Crop Sci. 48: 923-928.\\
Kindlmann, P., and A.F. Dixon. (2010). Modelling Population Dynamics of Aphids and Their Natural Enemies, pp. 1-20. In Kindlmann, P., A.F. Dixon, and J. Michaud (eds.), Aphid Biodiversity under Environmental Change. Springer, Dordrecht.\\
Kot, M. 2001. Elements of mathematical ecology. Cambridge University Press, UK.\\
Langa, J.A., J.C. Robinson, and A. Suárez, A. 2002. Stability, instability, and bifurcation phenomena in non-autonomous differential equations. Nonlinearity. 15: 887-903.\\
Liu, J., K.M. Wu, K.R. Hopper, and K.J. Zhao. 2004. Population dynamics of Aphis glycines (Homoptera: Aphididae) and its natural enemies in soybean in northern China. Ann. Entomol. Soc. Am. 97:235–39\\
O’Neal, M.E., A.J. Varenhorst, and M.C. Kaiser. 2018. Rapid evolution to host plant resistance by an invasive herbivore: soybean aphid (Aphis glycines) virulence in North America to aphid resistant cultivars. Curr. Opin. Insect. Sci. 26: 1-7.\\
Price, P.W., R.F. Denno, M.D. Eubanks, D.L. Finke, and I. Kaplan. 2011. Insect ecology: behavior, populations, and communities. Cambridge University Press, New York, NY.\\
Ragsdale, D.W., B.P. McCornack, R.C. Venette, B.D. Potter, I.V. MacRae, E.W. Hodgson, M.E. O’Neal, K.D. Johnson, R.J. O’Neil, C.D. DiFonzo, T.E. Hunt, P.A. Glogoza, and E.M. Cullen. 2007. Economic threshold for soybean aphid (Hemiptera: Aphididae). J. Econ. Entomol. 100: 1258-1267\\
Ragsdale, D.W., D.A. Landis, J. Brodeur, G.E. Heimpel, and N. Desneux. 2011. Ecology and management of the soybean aphid in North America. Annu. Rev. Entomol.  56: 375-399.\\
Ragsdale, D.W., D.J. Voegtlin, and R.J. O’Neil. 2004. Soybean aphid biology in North America. Ann. Entomol. Soc. Am. 97: 204-208.\\
Tabashnik, B.E., T. Brévault, and Y. Carrière. 2013. Insect resistance to Bt crops: lessons from the first billion acres. Nat. Biotechnol. 31: 510-521.\\
Tilmon, K.J., E.W. Hodgson, M.E. O'Neal, and D.W. Ragsdale. 2011. Biology of the soybean aphid, Aphis glycines (Hemiptera: Aphididae) in the United States. J. Integr. Pest Manag. 2: A1-A7.\\
Varenhorst, A.J., M.T. McCarville, and M.E. O’Neal. 2015. An induced susceptibility response in soybean promotes avirulent Aphis glycines (Hemiptera: Aphididae) populations on resistant soybean. Environ. Entomol. 44: 658-667.\\
Voegtlin, D.J., S.E. Halbert, and G. Qiao. 2004. A guide to separating Aphis glycines Matsumura and morphologically similar species that share its hosts. Ann. Entomol. Soc. Am. 97: 227–232. (2004).\\
Zhang, X. 1988. Study on the instars of soyabean aphid (Aphis glycines Matsumura) in Jilin Province. J. Jilin Agric. Univ. 10: 15-46.\\

\section{Tables and figures}

     \begin{figure}[!htbp]
        \includegraphics[width=\linewidth]{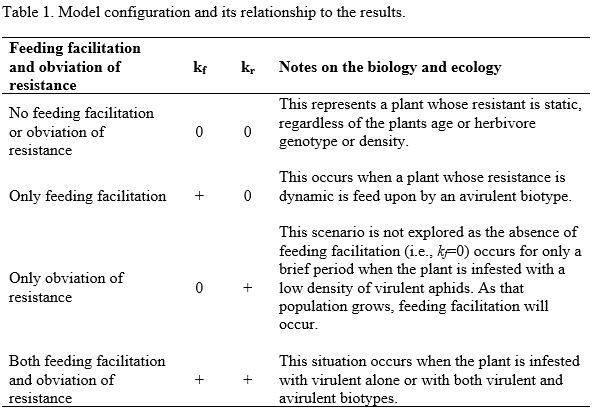} % first figure itself
        \label{fig:aubfig 1a}
        
    \end{figure}

         \begin{figure}[!htbp]
        \includegraphics[width=\linewidth]{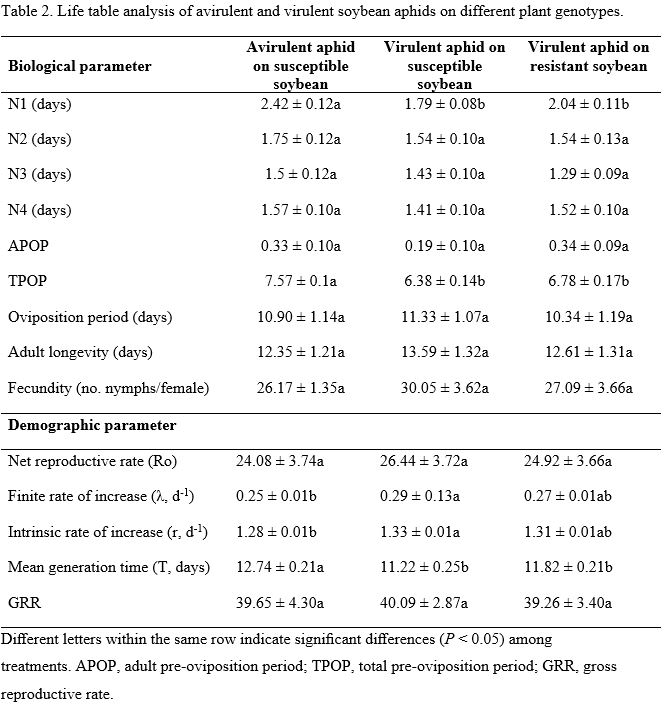} % first figure itself
        \label{fig:aubfig 1b}
        
    \end{figure}
    
     \begin{figure}[!htbp]
        \includegraphics[width=\linewidth]{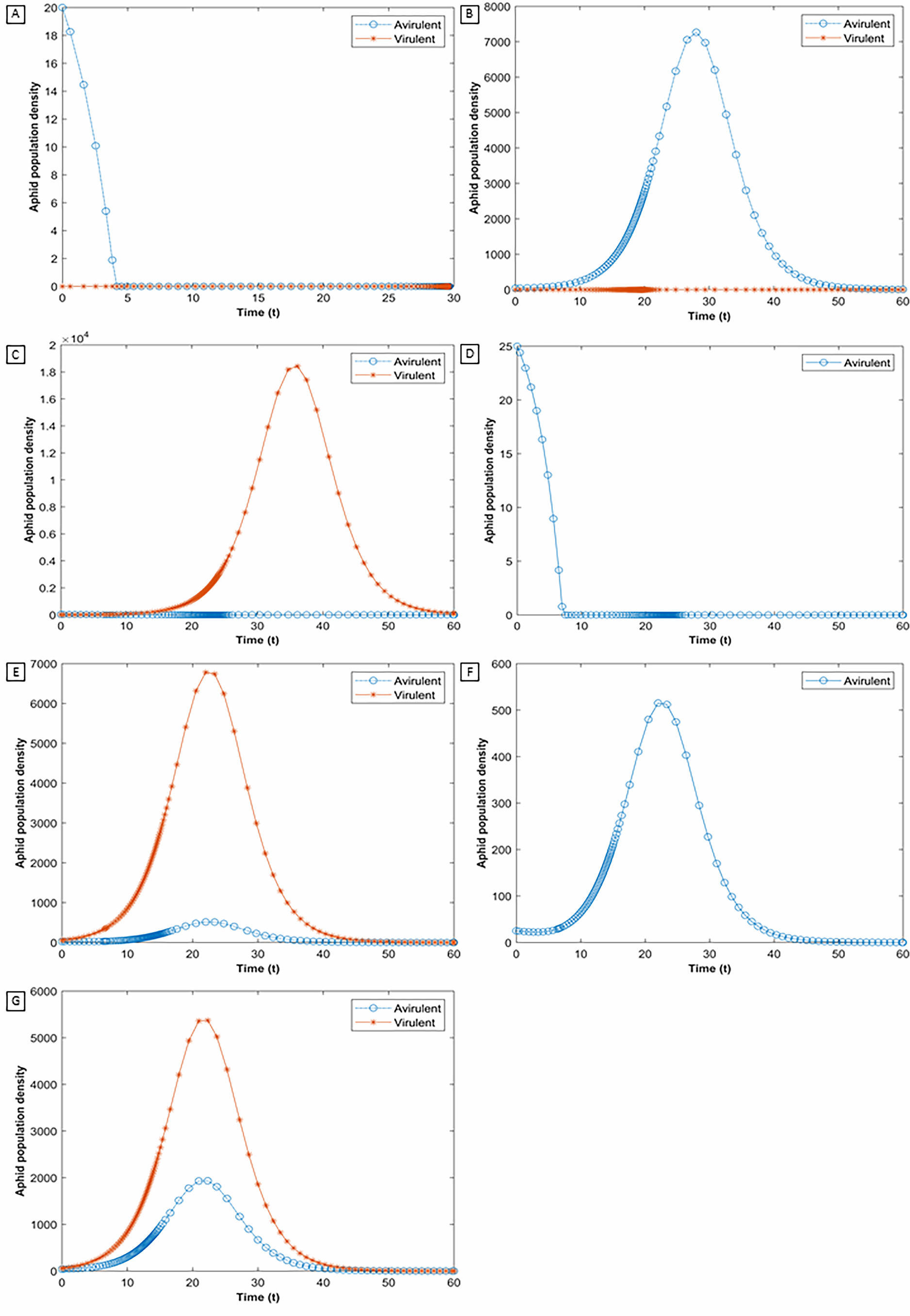} % first figure itself
        \caption*{{Figure 1. The impact of varying the initial population of the avirulent aphid ($x_A$), the virulent aphid ($x_V$) and the level of aphid-resistance in the host plant ($R$) are explored in 1A-1E. The initial populations of both biotypes are A) $x_A(0)$ =20 and $x_V (0)$ = 0; B) $x_A(0)$ =40 and $x_V (0)$ = 0; C and D) $x_A(0)$ =25 and $x_V (0)$ = 5; E and F) $x_A(0)$ =25 and $x_V (0)$ = 60 G) $x_A(0)$ =40 and $x_V (0)$ = 60. Results reported in 1A and 1C-1D demonstrate that the model accounts for the impact of resistance on an avirulent population, 1B and 1E-1G demonstrates the capacity for avirulent aphids to overcome this resistance. In 1E-1G, the virulent aphids reach a higher peak population than the avirulent aphids, with an increase in their initial population. In 1C and 1D avirulent aphid goes to extinction while virulent aphid reaches a high peak.}}
        \label{fig:aubfig 1a}
        
    \end{figure}
    
         \begin{figure}[!htbp]
        \includegraphics[width=\linewidth]{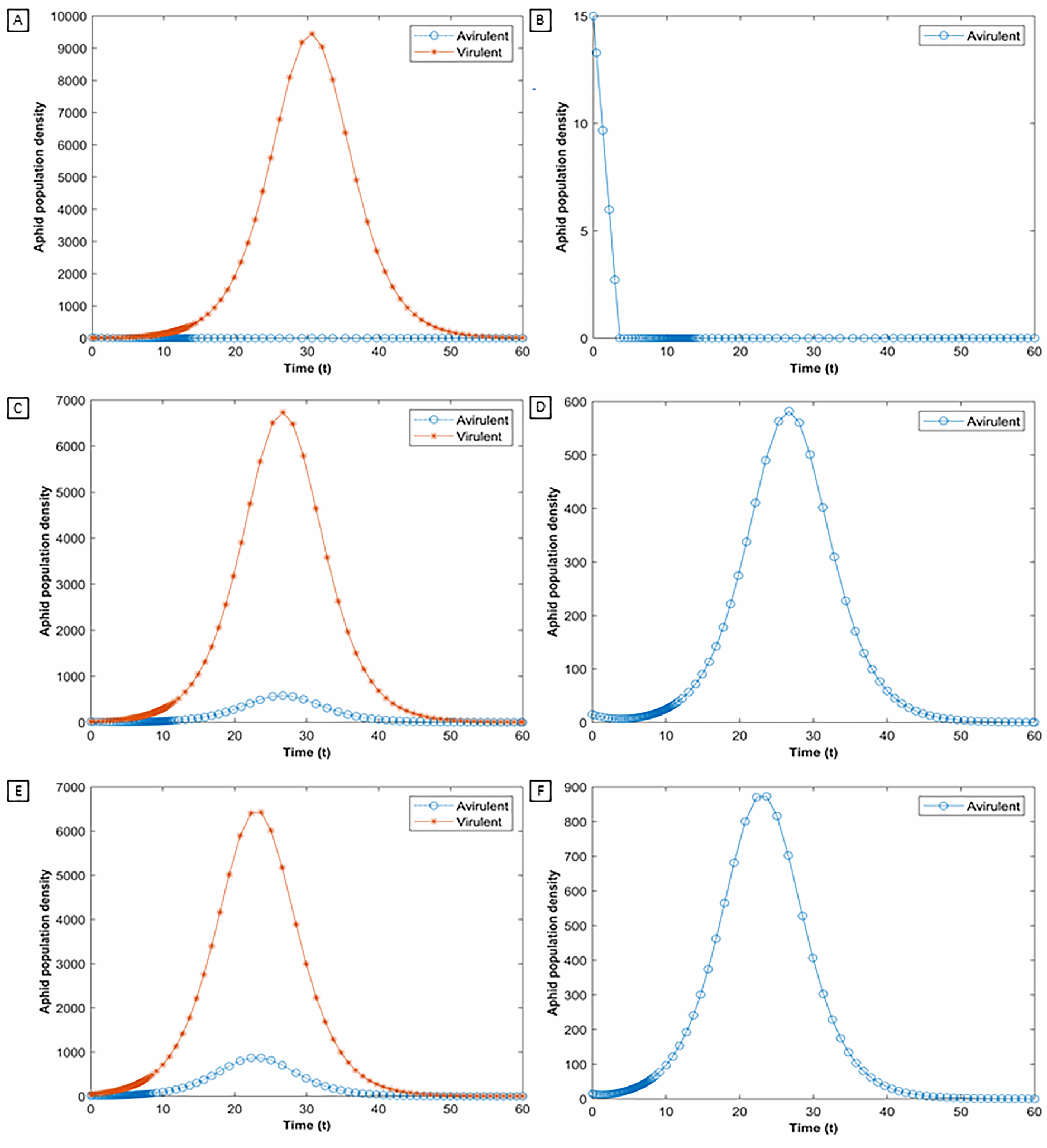} % first figure itself
        \caption*{{Figure 2. The impact of increasing the initial population of the virulent aphid (red line, $x_V$ ) on the co-existence of virulent and avirulent aphids (blue line, $x_A$) on a single, aphid-resistant soybean plant ($R(0)$ set to 30) are explored in 2A, 2C and 2E. The initial populations of both biotypes are A and B) $x_A(0)$ =15 and $x_V (0)$ = 10; C and D) $x_A(0)$ =15 and $x_V (0)$ = 20; E and F) $x_A(0)$ =15 and $x_V (0)$ = 50. In all the scenarios modeled, the virulent aphid out competes the avirulent aphid. However, as the initial population of virulent aphids increases, so too does the peak population of avirulent aphids. In figures 2B, 2D and 2F, the population of avirulent aphids is reported alone for each of the three scenarios, revealing their phenology across the modeled growing season. Note that the avirulent population went extinct early in the season when the initial population of virulent aphids was at its lowest (Fig. 2B).}}
        \label{fig:aubfig 1b}
        
    \end{figure}
     
         \begin{figure}[!htbp]
        \includegraphics[width=\linewidth]{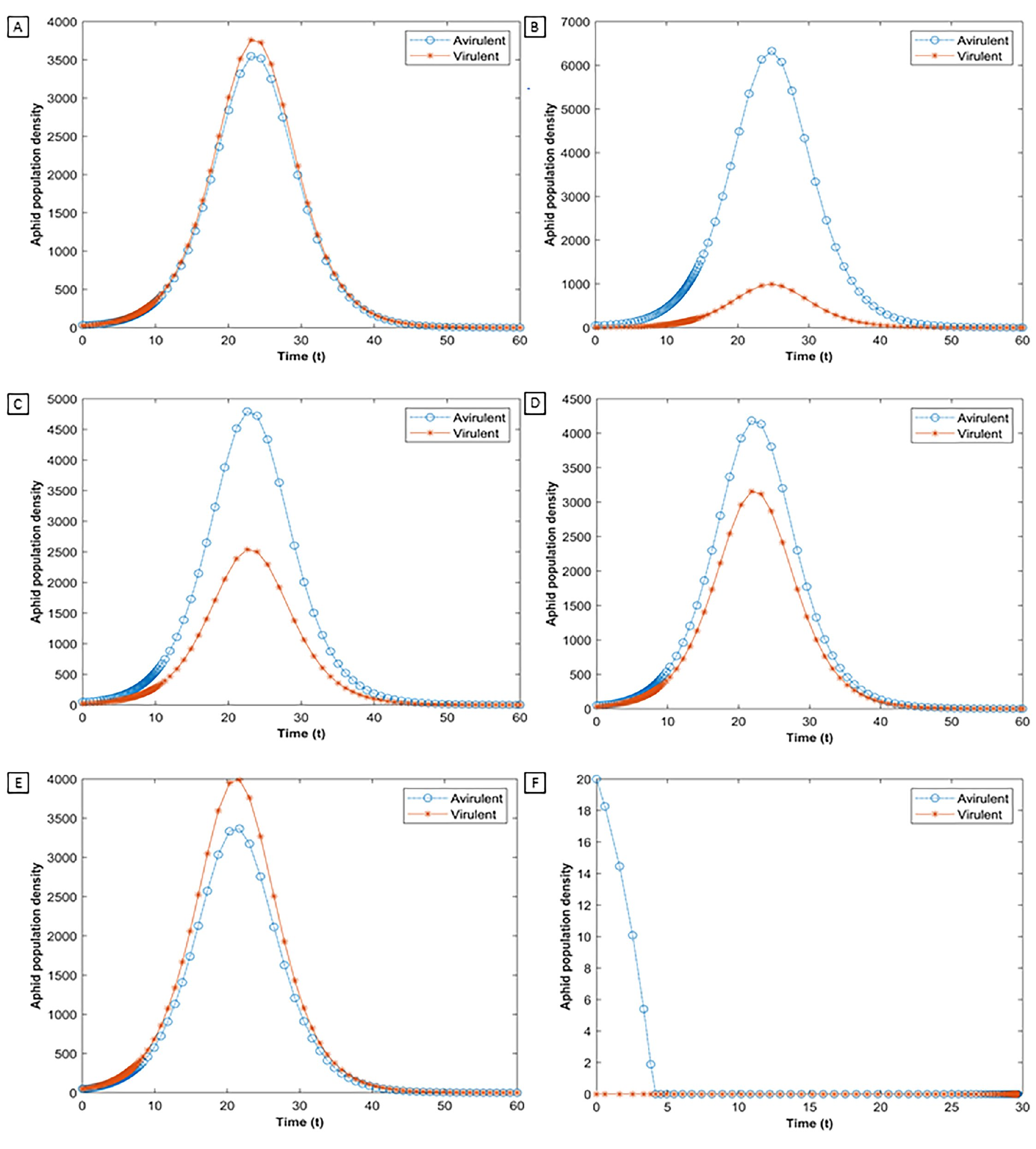} % first figure itself
        \caption*{{Figure 3. Soybean aphid dynamics when the initial populations of avirulent virulent aphid are greater than the resistance level. The initial populations of both biotypes are: A) $x_A(0)$ =35 and $x_V (0)$ = 25; B) $x_A(0)$ =50 and $x_V (0)$ = 5; C) $x_A(0)$ =50 and $x_V (0)$ = 20; D) $x_A(0)$ =50 and $x_V (0)$ = 30 E) $x_A(0)$ =50 and $x_V (0)$ = 50. }}
        \label{fig:aubfig 1c}
        
    \end{figure}

\end{document}